\documentclass{amsart}
\usepackage{amsmath, amsthm, amssymb, graphicx,slashed}
\def\cal{\mathcal}

\def\hat{\widehat}
\def\tilde{\widetilde}
\theoremstyle{definition}

\theoremstyle{remark}

\numberwithin{equation}{section} \font\teneurm=eurm10
\font\seveneurm=eurm7 \font\fiveeurm=eurm5
\newfam\eurmfam
\textfont\eurmfam=\teneurm \scriptfont\eurmfam=\seveneurm
\scriptscriptfont\eurmfam=\fiveeurm

 \font\teneusm=eusm10 \font\seveneusm=eusm7 \font\fiveeusm=eusm5
\newfam\eusmfam
\textfont\eusmfam=\teneusm \scriptfont\eusmfam=\seveneusm
\scriptscriptfont\eusmfam=\fiveeusm

\font\tencmmib=cmmib10 \skewchar\tencmmib='177
\font\sevencmmib=cmmib7 \skewchar\sevencmmib='177
\font\fivecmmib=cmmib5 \skewchar\fivecmmib='177
\newfam\cmmibfam
\textfont\cmmibfam=\tencmmib \scriptfont\cmmibfam=\sevencmmib
\scriptscriptfont\cmmibfam=\fivecmmib

\def\R{{\Bbb{R}}}

\begin{document}
\def\C{\Bbb{C}}
\def\M{{\cal M}}

\title{The Problem Of Gauge Theory}

\author{Edward Witten}
\address{Theory Group, CERN, Geneva Switzerland. On leave from
School of Natural Sciences, Institute for Advanced Study, Princeton
NJ 08540 USA.}
\email{witten@ias.edu}
\thanks{Supported in part by NSF Grant Phy-0503584.}


\date{November, 2008}

\def\Bbb{\mathbb}

\begin{abstract}
I sketch what it is supposed to mean to quantize gauge theory, and
how this can be made more concrete in perturbation theory and also
by starting with a finite-dimensional lattice approximation.  Based
on real experiments and computer simulations, quantum gauge theory
in four dimensions is believed to have a mass gap.  This is one of
the most fundamental facts that makes the Universe the way it is.
This article is the written form of a lecture presented at the
conference ``Geometric Analysis: Past and Future'' (Harvard
University, August 27-September 1, 2008), in honor of the 60th
birthday of S.-T. Yau.
\end{abstract}

\maketitle

\input epsf
\section{Yang-Mills Equations}
\label{intro}

\def\Tr{{\rm Tr}}

The strong, weak, and electromagnetic interactions -- in other
words, more or less everything we know about in nature except
gravity -- are all described by gauge theory.  Mathematically, in a
gauge theory with gauge group $G$, formulated on a spacetime $M$, a
gauge field is a connection on a $G$-bundle $E\to M$.  For our
purposes, $M$ is a four-manifold with a metric of Lorentz signature
$-+++$. In fact, for our purposes there is no essential loss to take
$M$ to be Minkowski spacetime $\Bbb{R}^{3,1}$ (that is, $ \Bbb{R}^4$
with a flat pseudo-Riemannian metric of Lorentz signature).

In addition to gauge fields, in nature there also are matter fields.
The matter fields describe things such as electrons, neutrinos,
quarks, and possibly Higgs particles.  Gauge fields mediate
``forces'' between particles described by matter fields (and between
additional particles described by the gauge fields themselves). For
simplicity, in this talk, I will omit the matter fields and just
describe the gauge fields.

One thing I should say before getting too far is that, for a large
variety of reasons, it is unrealistic, in a talk or a short article,
to expect to fully describe the problem of Yang-Mills theory. To
really appreciate the problem, it is necessary to delve into quantum
field theory in some depth. There is a large physics literature on
quantum field theory, and there is also a large math literature (for
example, see \cite{GJ}).

I felt in preparing the lecture that to make it comprehensible, I
could not simply skip the preliminaries.  One can skip the
preliminaries and make some formal statements, but such statements
are rather opaque.  On the other hand, it is also not possible to
fully explain the preliminaries in a short space or time, so I have
had to seek a middle path.  This inevitably involves some debateable
choices of what to explain.  The main idiosyncracy in my
presentation is that I have decided to rely on a Hamiltonian
approach (including the well-established \cite{Susskind} but
relatively unfamiliar notion of Hamiltonian lattice gauge theory)
rather than converting everything to path integrals. The path
integral approach is very powerful but involves an extra layer of
abstraction.

Classically, there is no problem to explain what is meant by
Yang-Mills gauge theory.  A gauge field, that is a connection $A$ on
a $G$-bundle $E\to M$, has a curvature $F=dA+A\wedge A$. The
curvature is a two-form on $M$ valued in the adjoint bundle ${\rm
ad}(E)$ derived from $E$. The classical Yang-Mills equations are
\begin{equation}\label{zurf} D\star F=0,\end{equation}
where $\star$ is the Hodge star (mapping two-forms to $n-2$-forms,
where $n={\rm dim}\,M$), and $D=d+A$ is the gauge-covariant
extension of the exterior derivative.

The Yang-Mills equations may be most familiar in the abelian case,
that is in the case that $G=U(1)$. Maxwell's equations of
electromagnetism (in vacuum) can be described in terms of a two-form
$F$ that obeys \begin{equation}\label{urf} 0 = d\star F=d
F.\end{equation}  The first of these equations is the Yang-Mills
equation as written in eqn. (\ref{zurf}). (Recall that $D$ reduces
to $d$ when $G$ is abelian.) We can omit the second equation if we
define $F=dA$ as the curvature of a connection $A$, for it is then
an identity -- the Bianchi identity.

For $G=U(1)$, Maxwell's equations describe propagation of light
waves in vacuum. These are linear equations, so, for example, in the
approximation that Maxwell's equations are valid, two beams of light
pass through each other without scattering.  (In a more precise
description of nature, there are all sorts of corrections to
Maxwell's equations, involving things such as quantum mechanics,
electrons, and gravity, so it is not expected that the propagation
of light is precisely linear even in vacuum.  The nonlinearities are
expected to be very small and have not yet been observed, though it
is believed that this may be possible in the near future.)

In the nonabelian case, the Yang-Mills equations $D\star F=0$ are
{\it nonlinear} wave equations.  So classical Yang-Mills waves do
scatter each other, although in the case of weak waves, the
nonlinearities and the scattering effects are small.  In this
respect, the Yang-Mills equations are analogous to the vacuum
Einstein equations $R_{\mu\nu}=0$ (where $R$ is the Ricci tensor).
They are nonlinear hyperbolic wave equations.

\begin{figure}
  \begin{center}
    \includegraphics[width=3in]{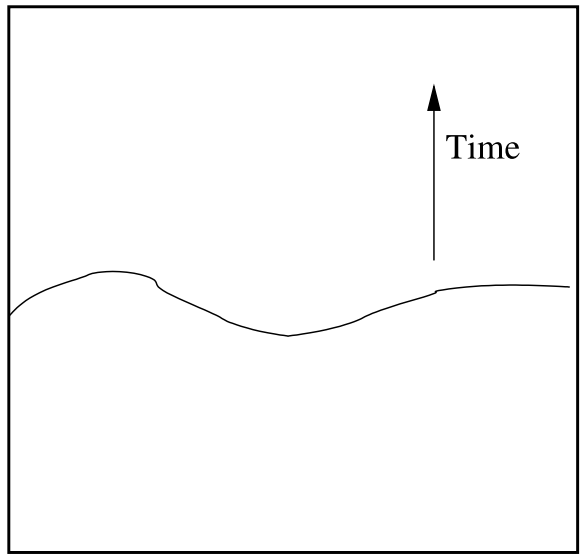}
  \end{center}
\caption{A spacetime with a global Cauchy hypersurface.}
  \label{Fig1}\end{figure}

Having a hyperbolic wave equation means that, optimistically
speaking, solutions of the equation are in one-to-one correspondence
with initial data, given on a global Cauchy hypersurface (fig.
\ref{Fig1}). However, there is a fundamental difference between
Yang-Mills theory and General Relativity. In General Relativity
there is a phenomenon of gravitational collapse -- the formation of
a black hole.  As a result, a solution on a Cauchy hypersurface
cannot be extended for all times, in general.     In Yang-Mills
theory, instead, there is a global-in-time existence result for
classical solutions (see \cite{Coleman}, \cite{CBS}).  This is
probably one of the reasons that quantum Yang-Mills theory is
simpler than quantum gravity, though it is not one of the reasons
that is easiest to interpret.

For gauge group $G=U(1)$, we observe classical solutions of
Maxwell's equations all the time -- light waves.  For nonabelian
$G$, even though there are beautiful nonlinear classical wave
equations, we do not observe these nonlinear classical waves in
practice.  That is actually because of a phenomenon known as the
mass gap.  The mass gap means that the description in terms of
nonlinear classical waves is only a good approximation above a
certain minimum energy and frequency.  At lower energies and
frequencies, one must use quantum field theory rather than classical
field theory.

According to theory, physical conditions that are well-described by
classical nonlinear wave equations can exist, but because of the
minimum frequency involved, our technology does not enable use to
generate the appropriate initial conditions. In practice, all
manifestations of Yang-Mills theory that we observe, except in the
abelian case, involve quantum behavior -- that is, they involve
phenomena that cannot be described by the classical field equations.
That is why the role of Yang-Mills theory in physics cannot be
described without talking about the quantum theory.

\section{Classical Phase Space}

\def\CW{{W}}
Formally speaking, the starting point in going to the quantum theory
is to observe that what I will call $\CW$, the space of all
solutions of the classical Yang-Mills equations modulo gauge
transformations, is an (infinite-dimensional) symplectic manifold.
The real reason for this is that the Yang-Mills equations are not
just equations.  They are the Euler-Lagrange equations associated
with an action function
\begin{equation}\label{folf} I=\frac{1}{4g^2}\int_M\,\Tr\,F\wedge
\star F. \end{equation} Here $\Tr$ is an invariant quadratic form on
the Lie algebra of $G$, and $g$ is a constant, known as the gauge
coupling constant.

In general, starting with any action, the space of classical
solutions of the Euler-Lagrange equations, modulo the relevant gauge
equivalence, is always a symplectic manifold.  Quantization has to
do with quantizing this symplectic manifold.

Actually, $\CW$ is a cotangent bundle.  This can be established as
follows.  Pick an initial value surface $S\subset M$, and let ${\cal
Y}$ be the space of all gauge fields on $S$ (that is all connections
on a $G$-bundle $E\to S$) modulo gauge equivalence.  Then $\CW$ is
the cotangent bundle of ${\cal Y}$, for any choice of $S$.  In
effect, to specify a classical solution of Yang-Mills theory
corresponding to a point in $\CW$, we must pick an initial value of
the gauge field along $S$ -- that is a point in ${\cal Y}$ -- and
also, as the Yang-Mills equations are of second order, we must
specify the normal derivative to the gauge field along $S$. By
forgetting the normal derivative, we get a map $\CW\to {\cal Y}$,
and $\CW$ is the cotangent bundle to ${\cal Y}$.

In finite dimensions, there is no problem in quantizing a cotangent
bundle.  But $\CW$ is infinite-dimensional, and in infinite
dimensions, we have to be careful.  For an elementary illustration
of the problem, recall that by a well-known theorem of Stone and von
Neumann, the quantization of $\Bbb{R}^{2n}$, with a symplectic
structure that comes from a nondegenerate skew form, is unique
(provided that one requires that this quantization should admit an
action of the Heisenberg group, the central extension of the group
of translations of $\Bbb{R}^{2n}$).  The analog of this for
$n=\infty$ is more subtle.

We need more information about what sort of answer we want to get.
The additional information is that the energy should be bounded
below.  With this information, quantization becomes unique again, at
least in the abelian case, as we will explain.

Instead of just describing the energy in an ad hoc way, let us
provide a framework for this discussion. Our Yang-Mills action
(\ref{folf}) is invariant under the symmetries of the
pseudo-Riemannian manifold $M$. Classically, we need only endow $M$
with a conformal structure (rather than a metric), since the action
is defined using the Hodge star operator $\star$, and in
four-dimensions, this operator is conformally-invariant when acting
on two-forms.  We take $M$ to be Minkowski spacetime $\R^{3,1}$,
with its standard conformal structure (induced from a flat
pseudo-Riemannian metric).  The group of conformal motions of $M$
(or more precisely of its conformal compactification) is $SO(2,4)$,
and this is a group of symmetries of the classical theory.

For $G=U(1)$, $SO(2,4)$ is realized as a group of symmetries of the
quantum theory, but this is actually not true for nonabelian $G$.
The quantum theory is obtained via a kind of limiting procedure,
concerning which I will try to give a few hints.  This limiting
procedure does not preserve the full $SO(2,4)$ symmetry, but only a
subgroup.  The details depend on exactly how one proceeds with
quantization.

\begin{figure}
  \begin{center}
    \includegraphics[width=2in]{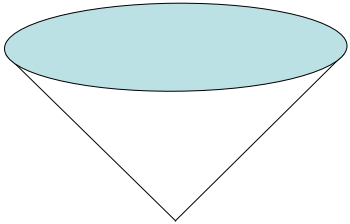}
  \end{center}
\caption{The future light cone.}
  \label{Fig2}\end{figure}

A maximal subgroup of $SO(2,4)$ that can be preserved in the
quantization, for nonabelian $G$, is the subgroup $P$ -- known as
the Poincar\'e group -- that preserves a flat metric on $\R^{3,1}$.
This group acts on linear coordinates $x$ on $\R^{3,1}$ by $x\to
ax+b$, where $b$ is a constant ``translation,'' and $a$ is a linear
transformation that belongs to the Lorentz group $SO(1,3)$ (the
group of symmetries of a quadratic form of Lorentz signature). We
are supposed to quantize in such a way that we get a unitary
representation of $P$.

Representation theory of $P$ is simple because $P$ is an extension:
\begin{equation}\label{dolf} 0\to \Bbb{R}^{3,1}\to P\to SO(3,1)\to
1.\end{equation}  Here $\Bbb{R}^{3,1}$ is the abelian group of
translations.  Its spectrum defines a point in the dual space to
$\Bbb{R}^{3,1}$.  This dual space, which I will denote
$\tilde\R^{3,1}$, is usually called momentum space.  Of course,
$\tilde\R^{3,1}$ is also endowed with a flat metric.  (We could use
the metrics to identify them, but this would be confusing.)

We write $p$ for a point in $\tilde\R^{3,1}$, usually called the
energy-momentum.  The condition $(p,p)=0$ defines a cone, called the
light cone.  $p$ is said to be lightlike if it lies on this cone. As
is usual in Lorentz signature, the light cone is the union of two
components, the ``past'' and ``future'' parts of the light cone. The
future light cone is sketched in fig. \ref{Fig2}.

We want a quantization such that the spectrum of energy-momentum
lies inside (and on) the future light cone.  This is usually
described by saying that the energy is bounded below by zero.  Here,
energy is a suitable linear function on $\tilde\R^{3,1}$.

Boundary points of the future light cone are allowed, but they play
a very special role.  The Hilbert space $\cal H$ that we get by
quantization is supposed to have a special state, the vacuum state
$|\Omega\rangle$, whose energy-momentum is supported at the apex of
the cone, in other words at the point $p=0$.  This state transforms
in a one-dimensional trivial representation of the Poincar\'e group.

\def\O{{\cal O}}
Apart from this one trivial representation, the other
representations that are relevant are constructed as follows.  Let
${\cal O}$ be a non-trivial $SO(3,1)$ orbit that is inside, or on,
the future light cone.    Let $V\to \O$ be a homogeneous vector
bundle over $\O$.  Then the space of ${\rm L}^2$ sections of $V$ is
a positive energy representation of the Poincar\'e group, and, apart
from the trivial representation that corresponds to the vacuum
state, these are the representations that we allow.

\begin{figure}
  \begin{center}
    \includegraphics[width=3.5in]{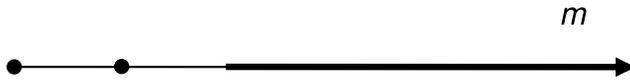}
  \end{center}
\caption{Sketched here is an example of what the spectrum of $m$
looks like when there is a ``mass gap.'' Apart from the isolated
eigenvalue at $m=0$ corresponding to the vacuum, there is in this
example a single discrete eigenvalue of $m$ at a positive value.
There is also a continuous spectrum (corresponding to multiparticle
states) which begins at twice the value of the smallest positive
eigenvalue.}
  \label{Fig3}\end{figure}

An orbit is characterized by the invariant $m^2=(p,p)$.  This
invariant is non-negative (since we consider only orbits inside or
on the light cone), and $m$, defined by the positive square root, is
called the mass. $m$ can have either a discrete or a continuous
spectrum. More exactly, in quantum field theory, $m$ always has a
continuous spectrum (from what are known as multiparticle states)
above some minimum value $m_*$. There  may also be a discrete
spectrum from single particle states with $m<m_*$. Positive discrete
eigenvalues of $m$ are called masses of particles.

The vacuum state always has $m=0$.  If this is a discrete eigenvalue
of $m$, in other words if there is $\epsilon>0$ such that every
state orthogonal to the vacuum has $m\geq \epsilon$, then we say
that the theory has a mass gap.  In fig. \ref{Fig3}, we sketch a
typical spectrum of a quantum field theory with a mass gap.  The
discrete eigenvalues of $m$ are 0 and one positive value.

Now, let us examine the problem of quantizing gauge theory, armed
with the information that the energy should be non-negative.
Classically, if we write the curvature in non-relativistic terms as
$F=dt\wedge d\vec x \cdot \vec E+\frac{1}{2}d\vec x\cdot d\vec
x\times \vec B$, where $\vec E$ and $\vec B$ are the electric and
magnetic fields, then the energy is the conserved quantity
\begin{equation}\label{zusk}H=\frac{1}{2g^2}\int
d^3x\,\Tr\,\left(\vec E^2+\vec B^2\right).\end{equation}  It is
non-negative, and vanishes precisely for the trivial solution with
$F=0$.

\section{Quantization}

We next discuss the quantization, beginning with the abelian case,
that is $G=U(1)$.  In the abelian case, the curvature $F$ is linear
in the connection $A$, that is $F=dA$, and Maxwell's equations
$d\star F=0$ are also linear.  The space $W$ of solutions modulo
gauge transformations is therefore also a linear space  -- an
infinite-dimensional one.

Now in the case of a {\it finite}-dimensional affine space,
$\Bbb{R}^{2n}$ for some $n$, quantization is unique (once one
requires that it should respect the affine structure, in a suitable
sense) according to a theorem of Stone and von Neumann.

This is far from being true in infinite dimensions.  Quantization of
an infinite-dimensional affine or linear space $W$ is far from
unique. But we do get the uniqueness again if we are given a
positive-definite quadratic function $Q$ on $W$, and we ask for a
quantization in which $Q$ is represented by a hermitian operator
that is bounded below.

In abelian gauge theory, because the curvature is a linear function,
the Hamiltonian or energy function $H$ is a quadratic function,
which moreover is positive-definite.   This puts us in the situation
just described.

\def\W{{\mathcal W}}
\def\H{{\mathcal H}}
The result of quantization can be described as follows.  As the
space of solutions of Maxwell's equations, $W$ has a natural
symplectic structure. This symplectic structure is
translation-invariant -- that is, it comes from a constant two-form
(on the infinite-dimensional linear space $W$). The
positive-definite quadratic function $Q$ on $W$ is equivalent to a
translation-invariant Riemannian metric on $W$. Combining these, $W$
is endowed with a translation-invariant complex structure and
therefore can be regarded as a complex vector space with a hermitian
metric, that is, a Hilbert space. Let us write $\W$ for $W$ regarded
in this way as a Hilbert space.

Quantization of $W$ is supposed to give us a Hilbert space $\H$.
This turns out to be the Hilbert space completion of the ``Fock
space'' constructed from $\W$, which by definition is
\begin{equation}\label{puik}\H_0=\oplus_{n=0}^\infty {\rm
Sym}^n\,\W.\end{equation} Here ${\rm Sym}^n\,W$ is the $n$-fold
symmetric product of $W$, with ${\rm Sym}^0\,W=\Bbb{C}$.

In (\ref{puik}), $\Bbb{C}={\rm Sym}^0\,W$ is the one-dimensional
space of ``vacuum'' states (Poincar\'e-invariant states of zero
energy, as discussed earlier).  $\W={\rm Sym}^1\,W$ is the space of
``single-particle states.''  Concretely, as a representation of the
conformal group, $\W$ is the space of sections of a certain
homogeneous vector bundle over the future light cone.  (In four
dimensions, this bundle is of rank 2, the number of polarization
states of an electromagnetic wave.)  Support on the cone means that
$\W$ is a space of massless states, that is states of $m^2=(p,p)=0$.
These are called the one-photon states.  Similarly, ${\rm Sym}^n\,W$
is the space of $n$-photon states.

This is our answer for quantization of abelian gauge theory, though
we need to say more about the physical interpretation in terms of
photons, and about how various classical expressions are realized as
operators acting on this Hilbert space.

\begin{figure}
  \begin{center}
    \includegraphics[width=2.5in]{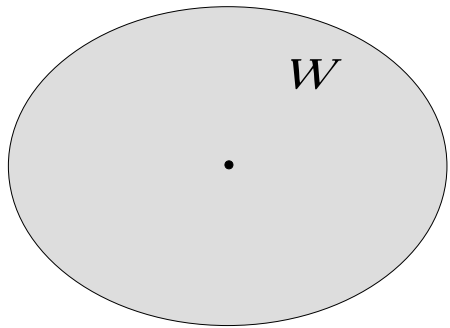}
  \end{center}
\caption{A space $W$ endowed with a function with a unique isolated
minimum, indicated by the black dot.}
  \label{Fig4}\end{figure}

Now let us discuss the nonabelian case.  The curvature $F$ is no
longer linear, so the space ${ W}$ of classical solutions is no
longer a linear or affine space.  Similarly, the energy function $H$
is not quadratic.  However, it is still true, as suggested in fig.
\ref{Fig4}, that $H$ is positive semidefinite with a unique zero
corresponding to the trivial solution $F=0$.

Like any function with an isolated minimum, $H$ looks quadratic near
its minimum.  One may ask whether this simple fact can be a starting
point for understanding the quantization.

The constant $1/g^2$ in $H$ is very important.  In general, suppose
that we have a not-quadratic function $H$ of variables $x_i$ with a
minimum at the origin:
\begin{equation}\label{luk}H=\frac{1}{g^2}\left(\sum_{i,j}a_{ij}x_ix_j
+\sum_{i,j,k}b_{ijk}x_ix_jx_k+\dots\right).\end{equation} We write
$y_i=x_i/g$, so that
\begin{equation}\label{uluk}H=\left(\sum_{i,j}a_{ij}y_iy_j
+g\sum_{i,j,k}b_{ijk}y_iy_jy_k+\dots\right).\end{equation}

We do not know how to diagonalize $H$ as an operator in a Hilbert
space, but we can diagonalize its quadratic part $H_0$:
\begin{equation}\label{luluk}H_0=\sum_{i,j}a_{ij}y_iy_j.\end{equation}
This means that if $g$ is small, we can approximately diagonalize
$H$.  The first step is to diagonalize $H_0$, and then one makes
successive corrections, treating the higher order terms in $H$ as
perturbations, so as to construct the eigenfunctions of $H$ in an
asymptotic expansion in powers of $g$.

Classically it does not make sense to say that $g$ is large or
small; $g$ is just an uninteresting constant multiplying the action.
But quantum mechanically there is a dimensionless number that in the
usual units is $g^2/\hbar c$.  This is really the small parameter in
the asymptotic expansion that was just suggested.  In this
asymptotic expansion, one diagonalizes $H$ -- and computes all
quantities of physical interest -- in an asymptotic expansion in
powers of $g^2/\hbar c$.

There is really a lot to explain here. There are many important
details and techniques in constructing the formal expansion, and
  there is actually much more to explain about what are the interesting and
  important things to calculate.  The techniques include Feynman
  diagrams, renormalization, path integrals, gauge fixing, and BRST
  symmetry.  What one wants to calculate are masses, other static
  quantities such as magnetic moments and other matrix elements of
  local operators, and especially scattering amplitudes.

  After a lot of work, one does end up with a systematic asymptotic
  expansion.  Moreover, there are many physics problems for which
  the asymptotic expansion is enough, in practice. That is the case
  for the electromagnetic and weak interactions, because $g^2/\hbar
  c$ is small (roughly 1/137 for electromagnetism).

Apart from the fact that the asymptotic expansion -- known as
perturbation theory -- is satisfactory for many questions, it is
important in another way: it is unrealistic to expect to develop an
exact theory without having a thorough understanding of perturbation
theory.  Trying to do this would be somewhat analogous to trying to
study Riemannian geometry without learning linear algebra.

\section{Nonperturbative Approach}

However, to understand the strong interactions, or nuclear force,
the asymptotic expansion is not enough.  We need to understand
something about the exact theory.  I will therefore conclude by
trying to say something about this.

Any known approach to understanding the exact theory requires, one
way or another, modifying it by introducing a ``cut-off'' so as to
make the number of variables effectively finite, and then taking a
limit in which the cutoff is removed. The usual way to do this is
via Euclidean lattice gauge theory, but we will not follow that
route because in this lecture we have avoided introducing path
integrals. (For an introduction to that subject, see for instance
\cite{GJ}.) Instead, I will describe the Hamiltonian version of
lattice gauge theory \cite{Susskind}.  This approach is not widely
used in practice, though possibly it could be useful.  At any rate,
whether or not this approach is useful in practice, it is easily
described and gives a good orientation about what it means to
introduce a cut-off and then remove it.

\def\A{\mathcal A}
\def\G{\mathcal G}
\def\Y{\mathcal Y}
A fact that was mentioned in section 2 is helpful here. This
concerns the space $W$ of solutions of the Yang-Mills equations
modulo gauge transformations. Pick an initial value hypersurface
$S$; we may as well simply take $S$ to be a ``time zero'' subspace
$\R^3\subset \R^4$. Let $\A$ be the space of gauge connections on
$S=\R^3$ and $\G$ the group of gauge transformations on $\R^3$. Thus
the quotient $\Y=\A/\G$ is the space of connections modulo gauge
transformations on $\R^3$. And $W$ can be identified as a cotangent
bundle, $W=T^*(\A/\G)$. The idea behind this identification is that
to determine a classical solution, we must give the initial value of
the connection (a point in $\A/\G$) and its time derivative (a
cotangent vector).   So to get a finite-dimensional approximation to
$W$, it suffices to get finite-dimensional approximations to $\A$
and $\G$.

\begin{figure}
  \begin{center}
    \includegraphics[width=2in]{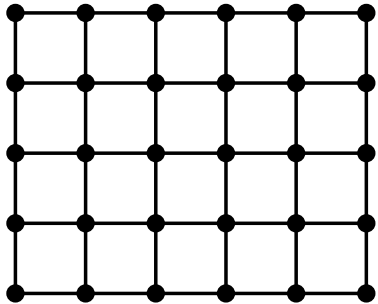}
  \end{center}
\caption{We introduce a cutoff by formulating gauge theory on a
lattice of finite spatial extent (for clarity, what is drawn here is
a two-dimensional rather than three-dimensional lattice). To
describe quantum field theory, we must take the spatial extent of
the lattice to infinity while also shrinking the lattice spacing to
zero.}
  \label{Fig5}\end{figure}

To do this, we approximate $\R^3$ by a finite set $\Gamma$ of points
-- later we will take the number of points to infinity.  In fact, as
in fig. \ref{Fig5}, we arrange the finite set of points in a regular
array, as part of a lattice.\footnote{\label{gong} Actually, it is
better to consider a discrete approximation to a three-torus, rather
than a discrete approximation to a parallelepiped in $\R^3$ as
sketched in the figure.  This avoids edge effects in the spectrum of
the Hamiltonian.  A discrete approximation to a three-torus is made
by taking an array of points that is finite and is periodic in three
directions.  The limit is still taken by letting the number of
points in the array go to infinity (so that it approximates a larger
and larger three-torus) while the distance between nearest neighbor
points goes to zero.} Of course, a finite set of points is not a
very good approximation to $\R^3$. To recover $\R^3$, we increase
the number of points.  In the process, we must improve the
approximation in two directions: we take the spacing between the
lattice points to be smaller and smaller, so as to recover the
continuum, while also taking the spatial extent of our chosen array
of points to be bigger and bigger, so that in the limit the array of
points covers all of $\R^3$.

So now, for each such finite set $\Gamma$, we must give an
approximation to $\A$, the space of all connections, and $\G$, the
group of gauge transformations.  Taking $\Gamma$ to be part of a
rectangular lattice, we connect the nearest neighbor pairs, as is
shown in the figure.  Then we approximate gauge theory by only
allowing parallel transport along lattice paths.  For each oriented
link $\ell$ between nearest neighbors in $\Gamma$, we introduce a
group element $U_\ell\in G$ that describes parallel transport along
$\ell$ (we take $U_{-\ell}=U_\ell^{-1}$).  By a connection on the
finite lattice corresponding to $\Gamma$, we mean a collection of
group elements $U_\ell$ for all $\ell$.  So our finite-dimensional
approximation to $\A$ is $\A_\Gamma=G^{n_1}$, where $n_1$ is the
number of nearest neighbor pairs in $\Gamma$.

Similarly, we permit ourselves to make gauge transformations at all
 points in $\Gamma$.  A gauge transformation is thus specified by
 giving an element $g_p\in G$ for all $p\in \Gamma$.  If $\ell$
 connects points $p$ and $q$, the gauge transformation acts on
 $U_\ell$ by $U_\ell\to g_pU_\ell g_q^{-1}$.  So our
 finite-dimensional approximation to the group $\G$ of gauge
 transformations is $\G_\Gamma=G^{n_0}$, where $n_0$ is the number
 of points in $\Gamma$.

\def\H{\mathcal H}
 The corresponding approximation to the space $W$ of classical
 solutions is $W_\Gamma=T^*(\Y_\Gamma)$, where
 $\Y_\Gamma=\A_\Gamma/\G_\Gamma$ is the space of gauge fields modulo
 gauge transformations.  Since we are in finite dimensions and
 $W_\Gamma$ is a cotangent bundle, quantization is straightforward:
 the Hilbert space associated to $\Gamma$ is just $\H_\Gamma={\rm
 L}^2(\Y_\Gamma)$.   Equivalently, this Hilbert space is the
 $\G_\Gamma$-invariant subspace of $\H^0_\Gamma={\rm L}^2(\A_\Gamma)$.
 Any operator on $\H^0_\Gamma$ that commutes with $\G_\Gamma$
 descends to an operator on $\H_\Gamma$.

 To complete the description of the finite-dimensional approximation, we
 must give a suitable approximation $H_\Gamma$ to the Hamiltonian or
 energy operator $H$.  It is not hard to do this.   It turns out
 that the first term in $H$, namely $H'=(1/2g^2)\int d^3x\,\Tr\,\vec
 E^2$, can be approximated by $(g^2/2)\Delta$, where $\Delta $ is
 the Laplace operator on the Riemannian manifold $\A_\Gamma$.
 ($\Delta$ commutes with $\G_\Gamma$, and so descends to an operator on
 $\H_\Gamma$ by virtue of the comment at the end of the last
 paragraph.)

 The second term in $H$, namely $H''=(1/2g^2)\int
 d^3x\,\Tr\,\vec B^2$, is the ${\rm L}^2$ norm of the curvature of a
 connection on the initial value surface.  To approximate this for
 the lattice $\Gamma$, the main point is to know what we mean by
 curvature in the context of such a finite-dimensional
 approximation.  This can be defined in terms of parallel transport around a small loop;
in the lattice $\Gamma$, the smallest possible nontrivial loops are
the squares of minimal area, which in lattice gauge theory are
usually called plaquettes.
 For any plaquette   $s\in \Gamma$, consisting of four nearest neighbor
 links $\ell_1,\dots,\ell_4$,
 we let $V_s$ be the function on $\A_\Gamma$ that
 associates to a connection $\{U_\ell\}$ the trace of the holonomy
 of that connection around $s$ (thus,
 $V_s=\Tr\,U_{\ell_1}U_{\ell_2}U_{\ell_3}U_{\ell_4}$; the trace is
 taken in the same representation used to define $H''$).
 Then $H''$ can be approximated by a suitable linear
 function of $\hat V=\sum_s V_s$.  (One subtracts a constant from
 $\hat V$ so that it vanishes when
 $U_\ell=1$ for all $\ell$; and one then multiplies by a constant to
 get an approximation to $H''$.)

We also want to define lattice approximations to other expressions
of classical and quantum gauge theory.  But I will not go farther;
the examples that have been given hopefully suffice to illustrate
the idea.

The problem of defining quantum gauge theory is, in this
formulation, to show that when we ``remove the cutoff'' by refining
and enlarging the finite set $\Gamma$, the lattice Hamiltonian
$H_\Gamma$ (and other operators) converge to a limit.\footnote{Here
to avoid edge effects, it is best to take $\Gamma$ to be an array of
points that is periodic in three directions, as mentioned in
footnote \ref{gong}.  Otherwise, there are edge effects at the
boundary of $\Gamma$ and the statement that $H_\Gamma$ converges to
a limit needs to be formulated carefully. } While removing the
cutoff, one must also adjust the coupling constant $g$ in a suitable
fashion.

There is a precise theory of how $g$ must be adjusted; if $a$ is the
lattice spacing (the distance between points in $\Gamma$), then one
requires $g\sim f/|\ln a|$ for $a\to 0$.  Here $f$ is a constant
that depends on $G$; it was computed by Gross, Wilczek, and Politzer
in 1973.  This computation led to the 2004 Nobel Prize for
``asymptotic freedom,'' which is the statement that $g$ must go to
zero as $a$ does.

\section{Breaking Of Conformal Invariance And The Mass Gap}

For  $a,g\to 0$, it is believed that there is a limiting theory that
obeys all of the axioms of quantum field theory, including
invariance under the Poincar\'e group.  However, it is believed that
in the limit one does {\it not} recover the $SO(2,4)$ conformal
symmetry of the classical theory.

A very basic aspect of the violation of the conformal symmetry is
that it is believed that the spectrum has a mass gap, and thus is
qualitatively as depicted in fig. \ref{Fig3}.  We recall that the
mass gap means simply that the mass $m$ of any state (orthogonal to
the vacuum) is bounded strictly above zero.

There is no mass gap in electromagnetism; the photon is massless, so
electromagnetic waves can have any positive frequency.  That is why
we can experience light waves in everyday life.  By contrast, the
mass gap in strong interactions means that the minimum frequency
needed to probe the world of $SU(3)$ gauge theory (which describes
the strong interactions) is $mc^2/\hbar$, where $m$ is the smallest
mass. Taking from experiment the value of the smallest mass, this
frequency is of order $10^{24}\,{\rm sec}^{-1}$, which is high
enough (with room to spare) that this world is way outside of our
ordinary experience.

 While it is very large compared to our ordinary experience,
 the mass gap is in one sense very small: it is zero in the asymptotic expansion
 described in Section 3.   As a result, we do not
 have a really good way to calculate it, though we know it is there
 from real experiments and computer simulations.

So in short, this mass gap is one of the most basic things that
makes the Universe the way it is, with electromagnetism obvious in
everyday life and other forces\footnote{The weak interactions are
also affected by a mass gap, but for very different reasons from the
strong interactions, which have been our subject here.  See \cite{W}
for an introduction.}  only accessible to study with modern
technology.

The mass gap is the reason, if you will, that we do not see
classical nonlinear Yang-Mills waves. They are a good approximation
only under inaccessible conditions.

I have spent most of my career wishing that we had a really good way
to quantitatively understand the mass gap in four-dimensional gauge
theory. I hope that this problem will be solved one day.

\bibliographystyle{amsplain}

\begin{thebibliography}{25}
\bibitem{GJ} J. Glimm and A. Jaffe, {\it Quantum Physics: A Functional
Integral Point Of View} (Springer, 1981).

\bibitem{Susskind} J. Kogut and  L. Susskind, ``Hamiltonian Formulation Of Wilson's Lattice
Gauge Theories,'' Phys. Rev. {
\bf D11} (1975) 395.

\bibitem{Coleman} S. Coleman,  ``There Are No Classical Glueballs,'' Commun. Math. Phys. {
\bf 55} (1977) 113.

\bibitem{CBS} Y. Choquet-Bruhat and I. Segal, ``Solution Globale Des \'Equations de Yang-Mills
sur l'Univers d'Einstein,'' C. R. Acad. Sci. Paris { \bf 294}
(1982) 225.

\bibitem{W}  E. Witten, ``From Superconductors And Four-Manifolds To Weak Interactions,"
Bull. Amer. Math. Soc.
{\bf 44}  (2007) 361.
\end{thebibliography}

\end{document}